\documentclass[12pt,twoside]{article}
\usepackage{amsmath,amssymb}

\title{On a fractional reaction-diffusion models arising in population dynamics} %
\author{S.H. Rasouli\\
 Department of Mathematics, Faculty of Basic Sciences,
\\Babol Noshirnani University of Technology, Babol, Iran\\e-mail:
s.h.rasouli@nit.ac.ir\\} 
\date{}

\begin{document}
\maketitle
\begin{center}
{\bf\large Abstract}\\
\end{center}

This paper is concerned with the existence of positive solutions for a
fractional population model with the homogeneous Dirichlet condition on the exterior of a
bounded domain. The approach is based on the sub-super solutions method. Our results generalize some recent results in the literature.\\\\
\hspace{-0.6 cm}Keywords: Fractional reaction-diffusion equation; Sub-supersolutions; Population models.
\section{Introduction}
\hspace{0.6 cm}In \cite{ekl-ss-rs1} the authors discussed the existence of multiple positive solutions to the steady state reaction diffusion model of the form
\begin{equation}
\left\{\begin{array}{ll}
-\Delta u  = \lambda \Big(u-\frac{u^{2}}{K}-c\frac{u^{2}}{1+u^{2}}\Big), & x\in \Omega,\\
u=0 ,  & x\in\partial \Omega,\\
\end{array}\right.
\end{equation}
where $\frac{1}{\lambda}$ is the diffusion coefficient, $K,$ and $c$  are positive constants and $\Omega\subset \mathbb{R}^{N}$ is a smooth bounded region with $\partial \Omega \in C^{2}.$ Here $u$ is the population density, $c$ is the maximum grazing rate. The term
$u-\frac{u^{2}}{K}$ represents a logistic growth, which means the per capita growth
 rate is linear depreciation. The term $\frac {u^{2}}{1+u^{2}}$ represents the rate of grazing by a constant number of grazing. Grazing type models apply to many ecological systems
arising in population dynamics such as the dynamics of salmon fish and spruce budworms.\\

 Recently in \cite{nf-rs-bs-ks}, authors extended the study of \cite{ekl-ss-rs1} to the following
 boundary value problem:
\begin{equation}
\left\{\begin{array}{ll}
-\Delta u  = \lambda \Big(u-\frac{u^{2}}{K}-c\frac{u^{2}}{1+u^{2}}\Big), & x\in \Omega,\\
\frac{\partial u}{\partial \eta}+ \sqrt{\lambda}u=0 ,  & x\in\partial \Omega,\\
\end{array}\right.
\end{equation}
where $\frac{\partial u}{\partial \eta}$ is the outward normal derivative of $u.$ The
authors established the occurrence of $S$-shaped bifurcation curves when parameters $c$ and $K$
satisfy certain conditions.\\

Most ecological systems also have some form of predation or harvesting of the population, which can either
be independent or dependent on the population density. Managers of ecosystems in which harvesting is well
regulated, such as fisheries management, often prefer constant yield harvesting (in which the harvesting rate
is density independent) over that of density dependent harvesting. In \cite{ekl-ss-rs2} the authors introduced a constant yield
harvesting term and studied the existence of positive solutions to the problem
\begin{equation}
\left\{\begin{array}{ll}
-\Delta u  = \lambda \Big(u-\frac{u^{2}}{K}-c\frac{u^{2}}{1+u^{2}}-\epsilon h(x)\Big), & x\in \Omega,\\
u=0 ,  & x\in\partial \Omega,\\
\end{array}\right.
\end{equation}
where the parameter $\epsilon \geq 0$ represent the level of harvesting, $h(x)\geq 0$
for $x\in \Omega,$ $h(x)=0$ for $x\in \partial\Omega,$ and $\|h\|_{\infty}=1.$ Here $\epsilon h(x)$ can be understood
as the rate of the harvesting distribution.\\

The goal of this paper is to extend
this study to the fractional population model of the form:

\begin{equation}
\left\{\begin{array}{ll}
(-\Delta)^{s} u  =  \lambda \Big(u-\frac{u^{2}}{K}-c\frac{u^{2}}{1+u^{2}}-\epsilon h(x)\Big), & \text{in} \,\,\Omega,\\
u =0 ,  & \text{in}\,\, \mathbb{R}^{N}\setminus\Omega,\\
\end{array}\right.
\end{equation}
where $-(\Delta)^{s} u$ is the fractional Laplacian operator of $u,$ $\Omega$ is a bounded domain in
$\mathbb{R}^{N}$ and $N>2s$ with $s\in (0,1),$  with sufficiently smooth boundary. Due
to this conditions, the extensions are challenging and nontrivial. One can refer to
\cite{gaa-shr, nf-jm-rs, shr1} for some existence results of population models.\\

Recently, a great deal of attention has been focused on studying of
problems involving fractional Sobolev spaces and corresponding
nonlocal equations, both from a pure mathematical point of view and
for concrete applications, since they naturally arise in many
different contexts, such as, among the others, the thin obstacle
problem, optimization, finance, phase transitions, stratified
materials, anomalous diffusion, crystal dislocation, soft thin
films, semipermeable membranes, flame propagation, conservation
laws, ultrarelativistic limits of quantum mechanics,
quasi-geostrophic flows, multiple scattering, minimal surfaces,
materials science and water waves. For more details, we can see
\cite{lac,ev,jlv}.\\

The natural space to look for solutions of the problem $(4)$ is the
usual fractional Sobolev space $W^{s,2}(\mathbb{R}^{N})=\Big\{u\in
L^{2}(\mathbb{R}^{N}): \Big(\frac{u(x)-u(y)}{|x-y|^{\frac{N}{2}+s}}\Big)\in
L^{2}(\mathbb{R}^{N} \times \mathbb{R}^{N})\Big\}$ endowed with the norm:
 \begin{equation}
 \|u\|_{W^{s,2}(\mathbb{R}^{N})} = \Big(\int_{\mathbb{R}^{N} \times \mathbb{R}^{N}}\frac{|u(x)-
 u(y)|^{2}}{|x-y|^{N+2s}}\,dx\,dy+\int_{\mathbb{R}^{N}}|u|^{2}\,dx\Big)^{1/2},
 \end{equation}
where the term
\begin{eqnarray*}
[u]_{W^{s,2}(\mathbb{R}^{N})} = \Big(\int_{\mathbb{R}^{N} \times \mathbb{R}^{N}}\frac{|u(x)-
 u(y)|^{2}}{|x-y|^{N+2s}}\,dx\,dy\Big)^{1/2}
\end{eqnarray*}
is the so-called Gagliardo (semi) norm of $u.$ To study
fractional Sobolev space in detail, we refer to \cite{ed-gp-ev,rs-ev}. We
define
\begin{eqnarray*}
X_{0}^{s,2}(\Omega)= \{u\in W^{s,2}(\mathbb{R}^{N}): u=0\,\,\text{a.e.}\, \text{in}\,\, \mathbb{R}^{N}\setminus \Omega\}.
\end{eqnarray*}
The space $X_{0}^{s,2}(\Omega)$ is a normed linear space endowed with
the norm $\|.\|_{X_{0}^{s,2}(\Omega)}$ defined as
\[
 \|u\|_{X_{0}^{s,2}(\Omega)} = \int_{\mathbb{R}^{N}\times \mathbb{R}^{N}}\frac{|u(x)-u(y)|^{2}}{|x-y|^{N+2s}}\,dx\,dy.
\]
We recall that, $X_{0}^{s,2}(\Omega)$ is a closed subspace of $W^{s,2}(\mathbb{R}^{N}),$ and its norm
is equivalent to the usual one defined in $(3).$\\

On the other hand, the spaces $W^{s,2}(\mathbb{R}^{N})$ and $X_{0}^{s,2}(\Omega)$  are strictly related to the fractional Laplacian operator. The fractional Laplacian is the pseudo-differential operator with Fourier symbol $\mathcal{F}$ satisfying
\begin{eqnarray*}
\mathcal{F} \Big((-\Delta)^{s} u\Big) (\xi) = |\xi|^{2s}\widehat{u}(\xi),\,\, 0<s<1,
\end{eqnarray*}
where $\widehat{u}$ denotes the Fourier transform of $u,$ (see \cite{ap}). Using Fourier transforms, it can be shown that (see \cite{ap}) an equivalent characterization of the fractional Laplacian is given by
\begin{eqnarray*}
(-\Delta)^{s} u(x)=C(N,s) P.V. \int_{\mathbb{R}^{N}}\frac{u(x)-u(y)}{|x-y|^{N+2s}}dy=C(N,s) \lim_{\epsilon \to 0^{+}}\int_{\mathbb{R}^{N}\setminus B_{\epsilon}(x)}\frac{u(x)-u(y)}{|x-y|^{N+2s}}dy.\\
\end{eqnarray*}
 Here $P.V.$ is a commonly used abbreviation for "in the principal value
sense"(as defined by the latter equation) and $C(N, s)$ is a dimensional
constant that depends on $n$ and $s,$ precisely given by
\begin{eqnarray*}
C(N,s)=\Big(\int_{\mathbb{R}^{N}}\frac{1-\cos \zeta_{1}}{|\zeta|^{N+2s}}\Big)^{-1},\,\,\zeta =(\zeta_{1},\zeta_{2},...,\zeta_{N})\in \mathbb{R}^{N}.\\
\end{eqnarray*}

In \cite{ed-gp-ev} (Proposition $3.6,$) the author proved the relation between the fractional Laplacian operator $(-\Delta)^{s}$and the Fractional Sobolev space $W^{s,2}(\mathbb{R}^{N}).$ They established
\begin{equation}
[u]_{W^{s,2}(\mathbb{R}^{N})}= 2 C(N,s)^{-1}\|(-\Delta)^{\frac{s}{2}}u\|_{L^{2}(\mathbb{R}^{N})}.
\end{equation}

\section{Main results}
\hspace{0.6 cm}In this section, we will establish our existence results. In order to state this existence result, we give some definitions. At first, we define a weak solution to the problem $(4)$ as follows.\\\\
{\bf Definition 2.1.} A function $u\in X_{0}^{s,2}(\Omega)$ is said to be a (weak) solution of $(4),$ if for any
 $\varphi \in X_{0}^{s,2}(\Omega),$ we have
\begin{equation}
\int_{\mathbb{R}^{N}}(-\Delta)^{\frac{s}{2}}u.(-\Delta)^{\frac{s}{2}}\varphi\,dx -\lambda \int_{\Omega} \Big(u-\frac{u^{2}}{K}-c\frac{u^{2}}{1+u^{2}}-\epsilon h(x)\Big)\varphi\, dx=0.\\
\end{equation}

Our main tool to prove the existence results is based on the application of sub- and supersolutions method. By a subsolution of $(4)$ we mean
a function $\underline{u}\in W^{s,2}(\mathbb{R}^{N})$ such that $\underline{u} \leq 0$ a.e. in $\mathbb{R}^{N}\setminus\Omega,$ and which
satisfies
\begin{equation}
\int_{\mathbb{R}^{N}}(-\Delta)^{s}\underline{u}\,.\,\varphi\,dx \leq \lambda \int_{\Omega} \Big(\underline{u}-\frac{\underline{u}^{2}}{K}-c\frac{\underline{u}^{2}}{1+\underline{u}^{2}}-\epsilon h(x)\Big)\varphi\, dx,
\end{equation}
for all  $\varphi \in X_{0}^{s,2}(\Omega)$ such that $\varphi \leq 0$ a.e. $\Omega.$ A supersolution $\overline{u}\in W^{s,2}(\mathbb{R}^{N})$
is defined analogously by reversing the inequality. Actually, it follows from \cite{mc-pg-eh} that if there exist a subsolution $\underline{u}$
and a supersolution $\overline{u}$ such that $\underline{u} \leq \overline{u}$ \,a.e. in $\Omega,$ then there exists a positive solution $u\in X_{0}^{s,2}(\Omega)$ satisfying $\underline{u} \leq u \leq \overline{u}$\, a.e. in $\Omega.$\\

Let $\lambda_{1}$ be the first eigenvalue of $(-\Delta)^{s}$ in $\Omega$ and $\phi_{1}>0$ be the corresponding
eigenfunction \cite{gmb-vdr-rs}, i.e.,

\begin{equation}
\left\{\begin{array}{ll}
(-\Delta)^{s} \phi_{1} =  \lambda_{1} \phi_{1}, & \text{in} \,\,\Omega,\\
\phi_{1} > 0,& \text{in} \,\,\Omega,\\
\phi_{1}=0,  & \text{in}\,\, \mathbb{R}^{N}\setminus\Omega.\\
\end{array}\right.
\end{equation}

The variational characterization of $\lambda_{1}$ is given by

\begin{eqnarray*}
\lambda_{1}= \inf_{v\in X_{0}^{s,2}(\Omega)\setminus \{0\}} \frac{\int_{\mathbb{R}^{N}}|(-\Delta)^{\frac{s}{2}}v|^{2}\,dx}{\int_{\Omega} v^{2}dx}.\\
\end{eqnarray*}

Moreover, it follows from \cite{xro} that there exist positive constants $c_{1},$ $c_{2}$ such that
\begin{equation}
0< c_{1}\delta^{s} \leq \phi_{1}(x)\leq c_{2}\delta^{s}\,\,\, \text{a.e}.\,\text{in}\,\, \Omega.\\
\end{equation}

We first note that if $\lambda \leq \lambda_{1},$ then $(4)$ has no positive solutions. This follows since if $u$ is a positive
solution of $(4),$ then $u$ satisfies
\begin{eqnarray*}
\int_{\mathbb{R}^{N}}|(-\Delta)^{\frac{s}{2}}u|^{2}\,dx\,=\,\lambda \int_{\Omega} \Big(u-\frac{u^{2}}{K}-c \frac{u^{2}}{1+u^{2}}-\epsilon h(x)\Big)u\,dx.
\end{eqnarray*}
But $\int_{\mathbb{R}^{N}}|(-\Delta)^{\frac{s}{2}}u|^{2}\,dx \geq \lambda_{1} \int_{\Omega} u^{2}dx.$ Combining, we obtain
\begin{eqnarray*}
\lambda \int_{\Omega} \Big(u-\frac{u^{2}}{K}-c \frac{u^{2}}{1+u^{2}}-\epsilon h(x)\Big)u\,dx \geq \lambda_{1} \int_{\Omega} u^{2}dx,
\end{eqnarray*}
 and hence
 \begin{eqnarray*}
(\lambda-\lambda_{1}) \int_{\Omega} u^{2}dx\geq \int_{\Omega} \Big(\frac{u^{2}}{K}+c \frac{u^{2}}{1+u^{2}}+\epsilon h(x)\Big)u\,dx \geq 0.
\end{eqnarray*}
This clearly requires $\lambda > \lambda_{1}.$\\

We now state our existence result.\\\\
{\bf Theorem 2.2.} Suppose that $\lambda > \lambda_{1}.$ Then there exists $\overline{\sigma}>0,$ $\underline{\sigma}>0,$ and $\epsilon^{*}>0$
such that for $K\in (\underline{\sigma},\overline{\sigma})$ and $\epsilon\in (0,\epsilon^{*})$ the problem $(4)$ has a positive solution.\\\\
{\bf Proof.} Let $\lambda > \lambda_{1}$ be fixed. We start with the construction of a positive supersolution. Let $\overline{u}=A e(x)$ where the constant $A>0$ is large
and to be chosen later and $e(x)$ is the unique positive solution of
\begin{equation}
\left\{\begin{array}{ll}
(-\Delta)^{s} e  = 1, & \text{in} \,\,\Omega,\\
e =0 ,  & \text{in}\,\, \mathbb{R}^{N}\setminus\Omega.\\
\end{array}\right.
\end{equation}
Moreover, it follows from \cite{gmb-vdr-rs} that there exist positive constants $\hat{c}_{1},$ $\hat{c}_{2}$ such that
\begin{equation}
0< \hat{c}_{1}\delta^{s} \leq e(x)\leq \hat{c}_{2}\delta^{s}\,\,\, \text{a.e}.\,\text{in}\,\, \Omega.
\end{equation}
We shall verify that $\overline{u}$ is a supersolution of $(4).$ To this end, let $\varphi \in X_{0}^{s,2}(\Omega)$ with $\varphi \geq 0.$ Since
$\max_{t\in \mathbb{R}}(t-\frac{t^{2}}{K})= \frac{K}{4},$ the inequality
\begin{eqnarray*}
Ae-\frac{(Ae)^{2}}{K}-c\frac{(Ae)^{2}}{1+(Ae)^{2}}-\epsilon h(x)\leq \frac{K}{4}
\end{eqnarray*}
hold for all $A\geq \frac{\lambda K}{4}$ a.e. in $\Omega.$ Then
\begin{eqnarray*}
\int_{\mathbb{R}^{N}}(-\Delta)^{s}\overline{u}\,.\,\varphi\,dx &=& A\,\int_{\Omega}\varphi\,dx\\
&\geq&\frac{\lambda K}{4} \int_{\Omega}\varphi\,dx\\
&\geq & \lambda \int_{\Omega}\Big(Ae-\frac{(Ae)^{2}}{K}-c\frac{(Ae)^{2}}{1+(Ae)^{2}}-\epsilon h(x)\Big)\varphi\,dx\\
&=& \lambda \int_{\Omega}\Big(\overline{u}-\frac{\overline{u}^{2}}{K}-c\frac{\overline{u}^{2}}{1+\overline{u}^{2}}-\epsilon h(x)\Big)\varphi\,dx.
\end{eqnarray*}
Hence, $\overline{u}=A e(x)$ is a supersolution of $(4).$\\

Next we construct a subsolution. Define $\alpha = \sqrt{\frac{\lambda_{1}}{\lambda}}\in (0,1).$ Then using $(10)$ and $(12),$ there exists $\theta>0$ such that $\phi_{1}(x)-\theta e(x)> \alpha \phi_{1}(x)>0$ a.e. in $\Omega.$ Let $\eta= 1+(\frac{1-\alpha}{2})^{2}.$ Define $$m_{\lambda}= \frac{1-\alpha}{2\|\phi_{1}(x)-\theta e(x)\|_{\infty}},$$
  $\overline{\sigma}=\frac{\eta}{c},$ $\underline{\sigma}=\frac{\eta}{2\eta+c}$ and $\epsilon^{*}=\frac{\theta m_{\lambda}}{\lambda}.$ We shall verify that $\underline{u}=m_{\lambda} (\phi_{1}(x)-\theta e(x))$ is a subsolution of $(4)$ for all $\epsilon\in (0,\epsilon^{*})$ and $K\in (\underline{\sigma},\overline{\sigma}).$ We note that $\lambda_{1}\leq \lambda \alpha\Big(1-\frac{(1-\alpha)}{2}(\frac{1}{K}-\frac{c}{\eta})\Big)$ follows immediately from the choice of $\alpha,$ $K$ and $\eta.$ Now, since $\theta m_{\lambda}-\lambda \epsilon^{*} =0$ by our choice of $\epsilon^{*},$ then $\theta m_{\lambda}-\lambda \epsilon =0$ for all $\epsilon <\epsilon^{*}.$
Let the test function $\varphi \in X_{0}^{s,2}(\Omega)$ with $\varphi \geq 0.$ Then it follows from $(9)$ that
\begin{eqnarray*}
&&\int_{\mathbb{R}^{N}}(-\Delta)^{s}\underline{u}\,.\,\varphi\,dx\\
 &=& m_{\lambda}\int_{\Omega}(\lambda_{1}\phi_{1}(x)-\theta)\,\varphi\,dx\\
 &\leq & \lambda \int_{\Omega}\Big[\alpha \phi_{1} m_{\lambda}(1-\frac{1-\alpha}{2\|\phi_{1}(x)-\theta e(x)\|_{\infty}})
 \|\phi_{1}(x)-\theta e(x)\|_{\infty}(\frac{1}{K}-\frac{c}{\eta})-\epsilon\Big]\,\varphi\,dx\\
 &\leq & \lambda \int_{\Omega}\Big((\phi_{1}-\theta e)m_{\lambda}\Big[1-m_{\lambda}(\phi_{1}-\theta e)(\frac{1}{K}-\frac{c}{1+(\frac{1-\alpha}{2})^{2}})-\epsilon\Big]\Big)\,\varphi\,dx\\
 &=&\lambda \int_{\Omega}\Big((\phi_{1}-\theta e)m_{\lambda}\Big[1-m_{\lambda}(\phi_{1}-\theta e)(\frac{1}{K}-\frac{c}{1+[m_{\lambda}\|\phi_{1}-\theta e\|_{\infty}]^{2}})-\epsilon\Big]\Big)\,\varphi\,dx\\
  &\leq & \lambda \int_{\Omega}\Big((\phi_{1}-\theta e)m_{\lambda}\Big[1-m_{\lambda}(\phi_{1}-\theta e)(\frac{1}{K}-\frac{c}{1+[m_{\lambda}(\phi_{1}-\theta e)]^{2}})-\epsilon\Big]\Big)\,\varphi\,dx\\
  &\leq & \lambda \int_{\Omega}\Big((\phi_{1}-\theta e)m_{\lambda}-\frac{[m_{\lambda}(\phi_{1}-\theta e)]^{2}}{K}
  -c\frac{[m_{\lambda}(\phi_{1}-\theta e)]^{2}}{1+[m_{\lambda}(\phi_{1}-\theta e)]^{2}}-\epsilon h(x)\Big)\,\varphi\,dx\\
   &=& \lambda \int_{\Omega}\Big(\underline{u}-\frac{\underline{u}^{2}}{K}
  -c\frac{\underline{u}^{2}}{1+\underline{u}^{2}}-\epsilon h(x)\Big)\,\varphi\,dx.
\end{eqnarray*}
Hence, $\underline{u}$ is a subsolution of $(4).$ Moreover, by choosing $A$ large we have $\underline{u} \leq \overline{u}.$ Hence there
exist a positive solution $u$ of $(4)$ such that $\underline{u} \leq u \leq \overline{u}$ and Theorem $2.2$ holds.\hspace{1 cm}$\Box$

\end{document}